\UseRawInputEncoding
\documentclass{amsart}
\usepackage{amsmath}
\usepackage{amssymb}
\newtheorem{thm}{Theorem}[section]
\newtheorem{cor}[thm]{Corollary}
\newtheorem{lem}[thm]{Lemma}
\newtheorem{prop}[thm]{Proposition}
\newtheorem{cons}[thm]{Construction}
\newtheorem{ex}[thm]{Example}
\theoremstyle{definition}
\newtheorem{defn}[thm]{Definition}
\newtheorem{notn}[thm]{Notation}
\theoremstyle{remark}
\newtheorem{rem}[thm]{Remark}
\newtheorem{fact}[thm]{Fact}
\long\def\Thm#1{\begin{thm} #1 \end{thm}}
\long\def\Cor#1{\begin{cor} #1 \end{cor}}
\long\def\Lem#1{\begin{lem} #1 \end{lem}}
\long\def\Prop#1{\begin{prop} #1 \end{prop}}
\long\def\Cons#1{\begin{cons} #1 \end{cons}}
\long\def\Def#1{\begin{defn} #1 \end{defn}}

\long\def\Rem#1{\begin{rem} #1 \end{rem}}
\long\def\Ex#1{\begin{ex} #1 \end{ex}}

\def\bar#1{\overline{#1}}
\def\Sect{\section}
\def\Rarr#1#2{\xrightarrow[#2]{#1}}

\long\def\Ref#1#2#3#4#5#6{
\bibitem{#1}
{\rm #2,}
\textit{#3.}
{\rm #4}
\textbf{#5}
{\rm #6.}
}
\long\def\Refb#1#2#3#4{
\bibitem{#1}
{\rm #2,}
\textit{#3.}
#4.
}
\long\def\comment#1{}
\def\Zz{{\mathbb Z}}
\def\Rr{{\mathbb R}}
\def\Cc{{\mathbb C}}
\def\Ff{{\mathbb F}}
\def\SS{\mathfrak{S}}

\def\ii{{\rm i}}
\def\e{{\rm e}}
\def\O{{\rm O}}
\def\into{\hookrightarrow}
\def\iso{\cong}
\def\leq{\leqslant}
\def\geq{\geqslant}

\def\comp{\mathbin{\mathchoice
{\circ}
{{\scriptstyle\circ}}
{{\scriptscriptstyle\circ}}
{{\scriptscriptstyle\circ}}
}}
\def\st{\mid}

\def\phi{\varphi}

\begin{document}

\title{Multivalued sections and self-maps of sphere bundles}
\author{M.~C.~Crabb}
\address{%
Institute of Mathematics,\\
University of Aberdeen, \\
Aberdeen AB24 3UE, UK}
\email{m.crabb@abdn.ac.uk}

\date{August 2022, revised May2024}
\begin{abstract}
Let $G$ be a finite group and $V$ a finite dimensional (non-zero)
orthogonal $G$-module such that, for each prime $p$
dividing the order of $G$, the subspace of $V$ fixed
by a Sylow $p$-subgroup of $G$ is non-zero and, if 
the dimension of $V$ is odd, has dimension greater
than $1$. Using ideas of Avvakumov, Karasev, Kudrya and Skopenkov
and work of Noakes on self-maps of sphere bundles,
we show that, for any principal $G$-bundle $P\to X$
over a compact ENR $X$, there exists a $G$-map
from $P$ to the unit sphere $S(V)$ in $V$.
\end{abstract}
\maketitle
\Sect{Introduction}
Let $G$ be a (non-trivial)
finite group and let $V$ be an 
orthogonal $G$-module of dimension $n\geq 1$. Given a principal
$G$-bundle $P\to X$ over a compact ENR $X$ we form the
$n$-dimensional real vector bundle $\xi =P\times_G V$
and its associated sphere bundle $S(\xi )=P\times_G S(V)$
(where $S(V)$ is the unit sphere in $V$).
\Thm{\label{gen}
Suppose that either
\par\noindent
{\rm (a)} for each prime $p$ dividing the order of $G$
the fixed subspace
$V^H$ of a Sylow $p$-subgroup $H$ of $G$ 
is non-zero and, if $n$ is odd,
has dimension greater than $1$, or
\par\noindent {\rm (b)} $n$ is odd, the fixed subspace
$V^H$ of a Sylow $2$-subgroup $H$ is non-zero
and $\dim X < 2(n-1)$.

Then there exists a $G$-equivariant map
$$
P \to S(V),
$$
or, equivalently, a section of the sphere bundle $S(\xi )$.
}
For a fixed integer $q>1$, the symmetric group $\SS_q$, the
group of permutations of the set $I=\{ 1,\ldots ,q\}$,
acts on the Euclidean vector space $\Rr [I]$ with 
orthonormal basis $e_i$, $i\in I$. 
Let $L$ denote the $(q-1)$-dimensional quotient of $\Rr [I]$
by the fixed subspace generated by $e_I=e_1+\ldots +e_q$.
We give $L$ the inner product for which the orthogonal
complement of $\Rr e_I$ projects isometrically onto $L$.
Given a prime $p$ we can write $q=\sum_ra_rp^r$, where $0\leq a_r<p$,
and set $\alpha_p(q)=\sum_r a_r$. Then, for a Sylow $p$-subgroup $H$
of $\SS_q$, the dimension of the fixed subspace $L^H$ is equal to
$\alpha_p(q)-1$. (See the discussion in Section \ref{sym}.)

Thus we deduce from Theorem \ref{gen}(a)
the following recent result of 
Avvakumov, Karasev, Kudrya and Skopenkov \cite{AK, AK1}, which was
applied in \cite{AKS} to the topological Tverberg theorem.
\Cor{\label{main}
{\rm (\cite{AK, AK1, AKS}.)}
Let $P\to X$ be a principal $\SS_q$-bundle over a compact ENR $X$
and let $\lambda$ denote the real vector bundle
$P\times_{\SS_q} L\to X$ of dimension $q-1$ over $X$.
Suppose that $q$ is not a prime power and that $k\geq 1$ is
a positive integer.
If $k=1$ and $q$ is even, suppose further that $\alpha_p(q)>2$
for every prime $p$.

Then there exists an $\SS_q$-equivariant map
$$
P \to S(\Rr^k\otimes L)
$$
to the unit sphere in $\Rr^k\otimes L$,
or, equivalently, a section of the sphere bundle
$S(\Rr^k\otimes\lambda )$ of $\Rr^k\otimes \lambda\,$.
\qed
}
The traditional approach to such problems, following
Sullivan \cite{mit}, is to localize at a
prime $p$.
If $H\leq G$ is a Sylow $p$-subgroup
(trivial if $p$ does not divide the order of $G$) 
and $V^H$ is non-zero,
then there is an $H$-equivariant map
$P\to S(V)$. Equivalently, the pullback $S(\pi^*\xi)$
to the finite cover $\pi : P/H\to X=P/G$ has a section.
Since the order of this cover, the index of $H$ in $G$,
is prime to $p$, it follows by a standard argument 
that the $p$-local obstruction to the existence
of a section of $S(\xi)$ is zero in the
metastable range: $\dim X < 2(n-1)$. If this holds
for all primes $p$,
the global obstruction is zero and the desired section,
or equivariant map, exists in this range of dimensions.
(An expository proof along these lines for the case of the symmetric
group $\SS_q$ was included in \cite[Proposition 6.5]{tver}.)
Integral cohomology is enough to establish the result
by the same reasoning in the range $\dim X \leq n$,
as in \cite[Theorem 5.1]{BZ2} and \cite[Theorem 1.2]{BZ}. 

The new proof in \cite{AK, AK1, AKS} without the 
dimensional restriction
is surprisingly elementary, especially when linked to a result of
Noakes \cite{lyle} on the degree of self-maps of sphere bundles.
This connection is explained in Section \ref{noakes}, where the main 
construction in \cite{AKS} is formulated in the language of
multivalued sections of sphere bundles and part (a) of 
Theorem \ref{gen} is established.
In Section \ref{sym} we discuss the extension (no longer a corollary)
of Corollary \ref{main} that was established in \cite{AK1} 
to the case in
which $k=1$ and $q$ is even but not twice a power of an odd
prime and give examples to show that this restriction
on $q$ is necessary.
The more technical Section \ref{stable}
describes the stable homotopy-theoretic
argument in the metastable range that was outlined above and 
completes the proof of part (b) of Theorem \ref{gen}. 
\Sect{\label{noakes}
Self-maps of sphere bundles}
We consider an orthogonal real vector bundle $\xi$ of dimension $n$ 
over a connected compact ENR $X$.
\Def{For a (fibrewise) 
self-map $f : S(\xi )\to S(\xi )$ of the 
sphere bundle of $\xi$ we define the
{\it degree of $f$}, $\deg (f)\in\Zz$,
to be the degree of the restriction of $f$ to any fibre,
which is independent of the fibre since $X$ is required
to be connected.}
We review first some results of Noakes \cite{lyle} from the 1970s, 
beginning with a basic construction.
\Def{\label{rho}
For a Euclidean vector space $V$ of dimension at least $2$
and an integer $r$, we define a continuous map
$$
\rho_r : S(V)\times S(V) \to S(V)
$$
as follows. Given $u,\ v\in S(V)$ lying in some $2$-dimensional
subspace $E\subseteq V$, we may choose a complex structure
on $E$ (compatible with the inner product)
and define $\rho_r(u,v)=z^ru$, where $v=zu$,
$z\in\Cc$, $|z|=1$. Since complex conjugation commutes with
the $r$th power, this
is independent of the choice of complex structure,
and, since $\rho_r(u,u)=u$ and $\rho_r(u,-u)=(-1)^ru$,
it is independent of the choice of $E$.
}
\Rem{More precisely, writing $\O (\Cc ,V)$ for the Stiefel
manifold of $\Rr$-linear isometric maps $\Cc\into V$, we
have a surjective map
$$
S(\Cc )\times \O (\Cc ,V) \to S(V)\times S(V),
\quad (z,a)\mapsto (a(1), a(z)),
$$
and $\rho_r$ lifts to the map
$$
\tilde\rho_r :
S(\Cc )\times \O (\Cc ,V) \to S(V),
\quad (z,a)\mapsto (a(1),a(z^r)).
$$
This provides a formal proof that $\rho_r$ is continuous.
}
When $V=\Cc$ we have the explicit formula
$$
S(\Cc )\times S(\Cc )\to S(\Cc ):
\rho_r(u,v)=u^{1-r}v^r,
$$
from which we see that, in general,
$\rho_r(v,u)=\rho_{1-r}(u,v)$ and, for a second integer
$s$, $\rho_{rs}(u,v)=\rho_r(u,\rho_s(u,v))$.
For any $V$, when $r=0$, the map $v\mapsto \rho_0(u,v)=u$ is constant;
and when $r=-1$, the map 
$v\mapsto -\rho_{-1}(u,v)=v-2\langle v,u\rangle u$ is
reflection in the hyperplane orthogonal to $u$.
\Lem{\label{degree}
Suppose that $\dim V$ is even.
Then, for a fixed $u\in S(V)$, the map $v\mapsto \rho_r(u,v):
S(V)\to S(V)$ has degree $r$.
}
The assertion is clear if $r=0$ and
when $\dim V=2$ it is evident from the explicit formula.
\begin{proof}
We give a proof from a
differentiable viewpoint \cite{milnor}.

Suppose that $\dim V >2$.
Let $E\subseteq V$ be the orthogonal complement of $\Rr u$.
Consider the smooth map $\pi : S(\Cc )\times S(E)\to S(V)$
taking $(x+y\ii ,e)$, where $x,\, y\in\Rr$, to $xu+ye$.
The restriction of $\pi : (S(\Cc )-\{ \pm 1\})\times S(E)
\to S(V)-\{ \pm u\}$ is a double cover, with covering involution
$(x+y\ii ,e)\mapsto (x-y\ii ,-e)$, which is orientation-preserving,
for an appropriate choice of orientations of
the connected manifolds $S(\Cc )\times S(E)$ and $S(V)$,
because $\dim V$ is even.

Now the map $f : v\mapsto \rho_r(u,v)$ lifts to the smooth map
$\tilde f:
(z,e)\mapsto (z^r,e): S(\Cc )\times S(E)\to S(\Cc )\times S(E)$.
Its degree is clearly $\deg \tilde f =r$.

The degree of $\pi$, for the chosen orientations,
is equal to $2$. Because the inverse image of a regular value $v$ in
the complement of $\{ \pm u\}$ consists of two points at each of which the sign of the determinant of $d\pi$ is $+1$.

Since $\deg (\pi )\deg (\tilde f)=\deg (\pi\comp \tilde f)
=\deg (f\comp \pi )=\deg (f)\deg (\pi)$, we see that $\deg f=r$.
\end{proof}
By the symmetry,
the degree of the map $u\mapsto \rho_r(u,v): S(V)\to S(V)$, 
for a fixed $v$, is equal to $1-r$.
\Prop{{\rm (Noakes \cite{lyle}).}
Let $\xi$ be an even-dimensional vector bundle. Then there is
a non-negative integer $N(\xi )$ such that there is a self-map
$S(\xi )\to S(\xi )$ of degree $1-d$ if and only if
$d$ is divisible by $N(\xi )$.
}
\begin{proof}
Consider two self-maps $f,\, g: S(\xi )\to S(\xi )$.
For each $r$ we may form $h=\rho_r(f,g): S(\xi )\to S(\xi )$.
Its degree, by Lemma \ref{degree}, is given by 
$$
\deg (h)=(1-r)\deg (f)+r\deg (g),
$$
that is, $1-\deg (h)= (1-r)(1-\deg (f))+r(1-\deg (g))$.
It follows that the set $I$ of integers of the form
$1-\deg (f)$ is an ideal.

\smallskip

{\small\par\noindent
(Indeed, $0\in I$ (because the identity map has degree $1$)
and $x,\,y\in I \Rightarrow (1-r)x+ry\in I$
for all $r\in\Zz$. If $I$ contains an element $y\not=0$, then 
$ry\in I$ for all $r\in\Zz$. So $I$ contains a least positive 
element $a$. Any element
of $I$ can then be written as $qa+b$ with $q\in\Zz$, $0< b\leq a$.
So $(1-r)(qa)+r(qa+b)=qa+rb\in I$ for all $r\in\Zz$. 
We must have $qa+sb-b< 0\leq qa+sb$ for some integer $s$, and
so $0\leq qa+sb<b\leq a$. Hence $qa+sb=0$ and $b=qa+(s+1)b\in I$.
By minimality, $b=a$.) 
}
\end{proof}
\Cor{Let $\xi$ be even-dimensional.
If $S(\xi )$ admits a section, then $N(\xi )=1$ and hence
$S(\xi )$ has self-maps of all degrees.
}
\begin{proof}
Given a section $s$, the self-map taking the constant value $s(x)$
in the fibre at $x\in X$ has degree $0$.
\end{proof}
\Prop{\label{properties}
{\rm (Some properties of $N(\xi )$).}
\par\noindent {\rm (i).}
Suppose that $\xi$ and $\eta$ are two even-dimensional
vector bundles over $X$. If $\xi$ can be embedded as a subbundle of
$\eta$, then $N(\eta )$  divides $N(\xi )$. 

\par\noindent {\rm (ii).}
The cohomology Euler class $e(\xi )\in H^n(X;\,\Zz (\xi ))$
with twisted integer coefficients is annihilated by $N(\xi )$,
that is,  $N(\xi )\cdot e(\xi )=0$.

\par\noindent {\rm (iii).}
Suppose that $\xi$ is a complex line bundle such that
$m\cdot c_1(\xi )=0\in H^2(X;\,\Zz )$ for some $m\geq 1$.
Then $N(\xi )$ divides $m$.
}
It follows from (ii) that, if the rational Euler class of $\xi$ 
is non-zero, then $N(\xi )=0$.
\begin{proof}
(i).
We can write 
$S(\eta )=S(\xi)*_X S(\xi^\perp)$ as the fibrewise join of
the sphere bundles of $\xi$ and its orthogonal complement 
$\xi^\perp$ in $\eta$. A fibrewise self-map $f$ of $S(\xi )$
gives a fibrewise self-map $f*_X1$
(the fibrewise join with the identity on $S(\xi^\perp)$)
of $S(\eta )$ with the same degree.

(ii). The homomorphism $H^n(X;\,\Zz (\xi ))\to H^n(X;\, \Zz (\xi ))$
induced by a self-map $f: S(\xi )\to S(\xi )$ of degree $d$ 
is multiplication
by $d$ and takes $e(\xi )$ to itself. So $(1-d)e(\xi )=0$.

(iii). Choose an isomorphism $\xi^{\otimes (m+1)} \to \xi$
between $\xi$ and its $(m+1)$-fold complex tensor product. 
The $(m+1)$th
power map $z\mapsto z^{\otimes (m+1)} : 
S(\xi )\to S(\xi^{\otimes (m+1)})\iso S(\xi )$ has degree $m+1$.
\end{proof}
\Lem{\label{sec}
Suppose that there exists a self-map $S(\xi )\to S(\xi )$
of degree $0$. Then the sphere bundle $S(\xi )$ admits
a section.
}
\begin{proof}
This is a standard nilpotence result.
Suppose that $f: E\to E$ is a self-map of a local trivially
fibre bundle over $X$ such that the restriction of $f$ to the
fibre $f_x : E_x\to E_x$ at $x\in X$ is homotopic to a constant
for every $x$. Then some power of $f$ is homotopic to a map
that is constant in each fibre.

\smallskip

{\small\par\noindent
(Using the local contractibilty of $X$ and compactness,
one gets a finite open cover $U_i$, $i=1,\ldots ,k$, of $X$
and maps $f_i$ homotopic to $f$ such that $f_i$ is constant
in the fibres over each point of $U_i$. Then the $k$th power
$f^k$ is
homotopic to the composition $f_1\comp \cdots\comp f_k$,
which is constant in every fibre.)  
}
\end{proof}
\Ex{Suppose that $\xi_1$ and $\xi_2$ are complex line bundles
over $X$ such that $m_1\cdot c_1(\xi_1)=0$ and 
$m_2\cdot c_1(\xi_2)=0$, where $m_1$ and $m_2$ are coprime
positive integers. Then $\xi =\xi_1\oplus\xi_2$ is isomorphic
to $\Cc\oplus (\xi_1\otimes\xi_2)$ as a complex bundle.
}
\begin{proof}
For $N(\xi )$ divides $N(\xi_i)$, which divides $m_i$,
for $i=1,\, 2$. Hence, $N(\xi )=1$ and so $S(\xi )$ has
a section. We conclude that $\xi$ splits as the direct
sum of a trivial line bundle $\Cc$ and a complex line
bundle classified by its first Chern class
$c_1(\xi )=c_1(\xi_1)+c_1(\xi_2)=c_1(\xi_1\otimes\xi_2)$.
\end{proof}
We look next at $d$-valued sections of a sphere bundle.
(A different class of multivalued sections is considered in
\cite[Section 5]{borsuk}.)
\Def{Let $\xi$ be an orthogonal real vector bundle over $X$.
For an integer $d\geq 1$, a {\it $d$-valued section}
of the sphere bundle $S(\xi )$ is a $d$-fold cover
$\tilde X\to X$ embedded, fibrewise over $X$, in $S(\xi )$.
}
Such a $d$-valued section determines, up to homotopy, a self-map of the
sphere bundle. To prepare for the proof, we recall a version of the
Pontryagin-Thom Umkehr construction introduced by McDuff in 
\cite{mcduff}.

Let $M$ be a connected, closed smooth manifold with 
tangent bundle $\tau M$.
We write $(\tau M)_M^+$ for the fibrewise one-point compactification
of $\tau M$; its fibre at $x\in M$ is the one-point compactification
$(\tau_xM)^+$ of the tangent space $\tau_xM$ at $x$.
It can be identified by stereographic projection with
the sphere bundle $S(\Rr\oplus\tau M)$ of the direct sum 
$\Rr\oplus\tau M$ of the tangent bundle and a trivial line bundle,
so that the point at infinity
in a fibre corresponds to the `North Pole' $(1,0)$:
$$\textstyle
v\in\tau_xM \mapsto 
(\frac{\| v\|^2}{4}+1)^{-1}(\frac{\| v\|^2}{4}-1,v)
\in S(\Rr\oplus\tau_xM).
$$

Consider a finite subset $N\subseteq M$ with cardinality $\# N=d$.
The restriction of $\tau M$ to $N$ will be denoted by $\nu$;
it is the normal bundle of the inclusion of the $0$-dimensional
manifold $N$ in $M$.
With a suitably scaled Riemannian metric on $M$ we have a tubular
neighbourhood $D(\nu )\into M$ embedding the closed unit disc bundle of the normal bundle $\nu$ of $N$ into $M$.
So $D(\nu)$ is the union of $d$ disjoint closed discs $D(\nu_y)$, 
$y\in N$. Over $D(\nu )\subseteq M$ we can identify the restriction 
of $\tau M$ with the pullback of $\nu$.
\Cons{\label{thom}
To the subset $N$, we associate a section $s_N$ of 
$(\tau M)^+_M$ 
taking the value $\infty$ (that is, $s_N(x)=\infty \in (\tau_xM)^+$)
outside the open unit disc bundle $B(\nu)$ and given on $D(\nu )$ by
$$
s_N(y,v)=[\psi (v)]\in \nu_y^+=\tau_{(y,v)}M^+\quad
\text{for $y\in N$, $v\in D(\nu_y)$,}
$$
where $\psi$ is the homeomorphism
$$
D(\nu_y)/S(\nu_y) \to \nu_y^+ :[v]\mapsto  v/(1-\| v\|^2)^{1/2}.
$$
}
We now look at the special case in
which $M$ is the sphere $S(V)$.
The tangent space $\tau_xS(V)$ at $x\in S(V)=M$ is the
orthogonal complement of $\Rr x$ in $V$.
The direct sum $\Rr\oplus\tau S(V)$ is thus identified with the trivial 
bundle $S(V)\times V$ (with $(1,0)\in\Rr\oplus\tau_xS(V)$ 
corresponding to $x\in V$), and so $(\tau M)_M^+=S(\Rr\oplus\tau M)$
is identified with the trivial sphere bundle $S(V)\times S(V)\to
S(V)$ (projecting onto the first factor)
in such a way that $\infty \in (\tau_xS(V))_{S(V)}^+$
corresponds to $(x,x)$.

The space of sections of the trivial bundle $S(V)\times S(V)\to S(V)$
is just the space of self-maps $S(V)\to S(V)$. And so
$s_N$ gives a map $f_N: S(V)\to S(V)$. In particular, if $N$ is empty
(and $d=0$), this gives the identity map $1: S(V)\to S(V)$.
\Lem{\label{deg}
The degree of $f_N$ is equal to $1-d$.
}
\begin{proof}
For a general connected manifold $M$ of dimension $m$, the homotopy
classes of sections $s$ of the sphere bundle $(\tau M)^+_M$ are 
classified by an integer $\deg (s)\in\Zz$.
To avoid twisted coefficients, let us assume that $M$ is oriented.
The sections $s$ are determinied by the difference class
$\deg (s)=\delta (s_\emptyset ,s)\in H^m(M;\,\Zz )=\Zz$, 
the obstruction to existence of a homotopy from the section
$s_\emptyset$ taking the value $\infty$ in each fibre to $s$.
For the section $s_N$, we have $\delta (s_\emptyset ,s_N)=d$
(as in \cite{mcduff}).

When we specialize to $M=S(V)$, the section $s_\emptyset$ corresponds
to the identity map $1=f_\emptyset :S(V)\to S(V)$ and
for any map $f: S(V)\to S(V)$ the difference class is equal to
$\delta (1 ,f)=1-\deg (f)=\deg (1)-\deg (f)\in\Zz$, 
which is the obstruction to
existence of a homotopy from $1$ to $f$.
Thus, $d=1-\deg (f_N)$ as required.

This argument is described in a more general setting in
Proposition \ref{mvs2}.
\end{proof}
The construction  that we have described is implicit in 
the arguments of \cite{AK, AK1,AKS}, although not there
linked to the, now classical, construction given in \cite{mcduff}.
\Prop{\label{mvs}
{\rm (Compare \cite[Lemma 10]{AKS}).}
Suppose that the sphere bundle $S(\xi )$
of an orthogonal real vector bundle $\xi$
over $X$ admits a $d$-valued section. 
Then there exists a self-map $S(\xi )\to S(\xi )$
of fibre degree $1-d$.
} 
\begin{proof}
Since $X$ is compact, for $\epsilon >0$ sufficiently small,
the closed balls $D_\epsilon (\tilde x)$ of radius 
$\epsilon$ (in the Euclidean metric) centred at the points
$\tilde x\in \tilde X$ in any fibre are disjoint.
So the construction of $f_N$ for a finite set $N$ can be carried
through fibrewise over $X$ to produce a map $f_{\tilde X}:
S(\xi )\to S(\xi )$ with fibre degree $1-d$.
\end{proof}
\Rem{When $d=2$ we can write down explicitly
a map $f:S(\xi )\to S(\xi )$
of degree $1-2=-1$ specified in each fibre by the 
(linear) reflection 
that interchanges the two points of $\tilde X$ in the fibre.
}
\Rem{When $\xi$ is a complex line bundle, 
the cover $\tilde X\to X$ determines a section of 
the sphere bundle $S(\xi^{\otimes d})$ of the complex
$d$-fold tensor power by taking the product of the $d$ elements
in the fibre of $\xi$. The section trivializes the line bundle
$\xi^{\otimes d}$ and so gives a self-map of degree
$1-d$ as in Proposition \ref{properties}(iii).
}
This is already enough to establish Theorem \ref{gen}(a) when $n$ is even.
\begin{proof}[Proof of Theorem \ref{gen}(a) when $n$ is even]
For a prime $p$ dividing $\# G$, choose a point $x\in S(V)$ such that
the stabilizer $G_x$ of $x$ contains a Sylow $p$-subgroup
of $G$. The cardinality $d$ of the orbit $Gx$ is, thus, not divisible
by $p$. Then $\tilde X =P\times_G Gx\subseteq P\times_G S(V)=S(\xi )$
gives a $d$-valued section of $S(\xi )$.
By Proposition \ref{mvs}, $S(\xi )$ has a self-map of fibre degree
$1-d$. Hence $N(\xi )$ divides $d$, which divides $\# G$. 
Since $p$ does not divide $d$, it does not divide $N(\xi )$. 

Hence $N(\xi )=1$, because $N(\xi )$ divides $\# G$
and is not divisible by any prime $p$ dividing $\# G$.
Thus $S(\xi )$ has a self-map of fibre degree $0$ and
so has a section, by Lemma \ref{sec}.
\end{proof}
\Prop{\label{dickson}
Suppose that $G$ is an elementary abelian
$p$-group of order $p^s$, $s\geq 1$, 
and $V$ is a $G$-module such that $V^G=0$. Then there
exists a principal $G$-bundle $P\to X$ admitting no
$G$-map $P\to S(V)$.
}
\begin{proof}
It is enough to consider a representation 
$V =\Cc^k\otimes (\Rr [G]/\Rr )$, where $k\geq 1$.
Take $G=C_p\times \cdots \times C_p$, where $C_p$ is
the group of $p$th roots of unity in $\Cc$
acting on $P= S(\Cc^m)\times \cdots \times S(\Cc^m)$
by $(g_i)\cdot (v_i)=(g_iv_i)$. For $m$ sufficiently large
the $\Ff_p$-Euler class $e(\xi )\in H^n(X;\,\Ff_p)$
of the complex vector bundle $\xi$ is non-zero.
\end{proof}
\Rem{An argument using complex $K$-theory Euler classes,
due to Munkholm \cite{munkholm}, shows similarly that,
if $G$ is cyclic of prime power order $p^s$ 
and $V$ is a $G$-module such that $V^G=0$, then there
exists a principal $G$-bundle $P\to X$ admitting no
$G$-map $P\to S(V)$.
The same result for any $p$-group $G$
\footnote{I am grateful to Roman Karasev for pointing this out.}
can be established by
using the Euler class in stable cohomotopy
and Carlsson's verification of Segal's Burnside ring conjecture
\cite{carlsson}.
(In the notation of Section \ref{stable}, the equivariant
stable cohomotopy class $\gamma (V)\in\omega^0_G(*;\, -V)$
of $V$ generates a free summand $\Zz \gamma (V)$, because its
image under the $G$-fixed-point map $\omega^0_G(*;\, -V)
\to \omega^0(*)=\Zz$ is $1$. By the Segal conjecture, there is
a bundle $P\to X$ such that the Euler class
$\gamma (\xi )\in\omega^0(X;\, -\xi )$, corresponding to the image of $\gamma (V)$ in $\omega^0_G(P;\, -V)$,
is non-zero; the sphere bundle $S(\xi)\to X$ has no section.)
}
To complete the proof of Theorem \ref{gen}(a), we need
a generalization of Proposition \ref{mvs}.
We start by extending Construction \ref{thom}.
\Cons{\label{thom1}
Suppose that we are given a family $h$ of pointed
maps $h_x : \nu_x^+\to \nu_x^+$ for $x\in N$.
To the subset $N$ and family $h$, we associate a section $s_{N,h}$
of $(\tau M)^+_M$ taking the value $\infty$
outside the open unit disc bundle $B(\nu)$ and given on $D(\nu )$ by
$$
s_{N,h}(y,v)=[h_y (\psi (v))]\in \nu_y^+=\tau_{(y,v)}M^+\quad
\text{for $y\in N$, $v\in D(\nu_y)$.}
$$
}
When $M=S(V)$ is a sphere, the section $s_{N,h}$ corresponds
to a self-map $f_{N,h} : S(V)\to S(V)$, which 
(by the argument used to prove Lemma \ref{deg}) has
degree $1-\sum_{y\in N} \deg h_y$.

Suppose that $\pi : \tilde X\to X$ is a finite cover defining
a $d$-valued section of $S(\xi)$. Then the fibre of $\pi^*\xi$
at $u\in\tilde X$ splits as an orthogonal direct sum $\Rr u\oplus
\eta_x$ and $\pi^*\xi$ is, thus, a direct sum $\Rr\oplus\eta$ of
a trivial line bundle and a real vector bundle $\eta$ of
dimension $n-1$.
Suppose that the space $\tilde X$ is expressed
as a disjoint union of subspaces
$\tilde X_\gamma$, $\gamma\in\Gamma$ such that the restriction 
$\tilde X_\gamma \to X$ of $\pi$
is a finite cover, of order $d_\gamma$ say.
(For example, the $\tilde X_\gamma$ might be the connected components
of $\tilde X$.)
\Prop{\label{mvs1}
{\rm (Compare \cite[Lemma 3.1]{AK1}).}
In the notation introduced above, suppose that $h : S(\eta )\to 
S(\eta)$ is a self-map of the sphere bundle $S(\eta )$ over $\tilde X$
such that, for each $\gamma\in\Gamma$,
the restriction of $h$ to $\tilde X_\gamma$
has constant fibre degree, $m_\gamma$ say.

Then $S(\xi )$ admits a self-map of fibre degree
$$
1-\sum_{\gamma\in\Gamma} m_\gamma d_\gamma \, .
$$
}
\begin{proof}
The bundle map $h: S(\eta )\to S(\eta )$ extends radially to a
map $\eta^+_{\tilde X} \to \eta^+_{\tilde X}$.
Then the construction of $f_{N,h}$ above can be performed fibrewise
over $X$ to produce the required self-map $f_{\tilde X,h}$
of $S(\xi )$ as in Proposition \ref{mvs}.
\end{proof}
\begin{proof}[Proof of Theorem \ref{gen}(a) when $n$ is odd]
Consider a prime $p$ dividing $\# G$ and a Sylow $p$-subgroup $H$
of $G$. Since $S(V^H)$ is an infinite set and $G$ has only finitely 
many subgroups, there exist two points $x,\, y\in S(V^H)$ such
that $G_x=G_y$ and $y\notin \{ \pm x\}$. It is possible that
$G_x$ contains a Sylow $p'$-subgroup for some other prime $p'$.

We can, therefore, choose finitely many points $x_\gamma ,\,
y_\gamma$, indexed by $\gamma\in\Gamma$, such that:
(i) $x_\gamma\not=
\pm y_\gamma$ and $G_{x_\gamma}=G_{y_\gamma}$;
(ii) for each
prime $p$ dividing $\# G$ at least one of the subgroups
$G_{x_\gamma}$ contains a Sylow $p$-subgroup of $G$;
(iii) if $\gamma\not=\delta$, then the subgroups $G_{x_\gamma}$
and $G_{x_\delta}$ of $G$ are not conjugate.

Let $\tilde X_\gamma =P\times_G Gx_\gamma \subseteq 
P\times_GS(V)=S(\xi)$.
Then $\tilde X_\gamma \to X$ is a finite cover of order
$d_\gamma =\# Gx_\gamma$.
If $\gamma\not=\delta$, the sets $\tilde X_\gamma$ and 
$\tilde X_\delta$ are disjoint, because the stabilizers
$G_{x_\gamma}$ and 
$G_{x_\delta}$ are not conjugate. Take 
$\tilde X=\bigsqcup_{\gamma\in\Gamma} \tilde X_\gamma$.
By construction, $d_\gamma$ divides $\# G$ and for each prime
$p$ dividing $\# G$ there is some $\gamma$ such that $p$ does
not divide $d_\gamma$.
So there are integers $m_\gamma$ such that 
$1-\sum_\gamma m_\gamma d_\gamma =0$.

Now $y_\gamma$ gives a section of $S(\pi^*\xi )=S(\Rr\oplus\eta)$
over $\tilde X_\gamma$: $[a,gx_\gamma ]\mapsto [a,gy_\gamma]$
for $a\in P,\, g\in G$
(well-defined, because $G_{x_\gamma}=G_{y_\gamma}$).
Then, because $gy_\gamma\not=\pm gx_\gamma$, the $\eta$-component
is non-zero and so determines a section of $S(\eta)$ over 
$\tilde X_\gamma$.
Hence, there is a map $S(\eta )\to S(\eta )$ with degree 
$m_\gamma$ over $\tilde X_\gamma$.
\end{proof}
It is easy to write down a counterexample if the assumption that
the dimension of $V^H$ is greater than $1$ is omitted. 
\Prop{\label{counter}
Let $G$ be cyclic of order
$6$, generated by an element $g$, and let $A$ and $B$ be the
$G$-modules $A=\Rr$, $B=\Cc$
with $g$ acting as $t\mapsto -t$, $z\mapsto\e^{2\pi \ii /3} z$,
respectively. Take $V=A\oplus B$ and
let $P$ be the $6$-dimensional free $G$-space
$S(\Rr^2\otimes A)\times S(\Rr^3\otimes B)$. 
Then

{\rm (i).}
There is a $G$-map $S(V)\to S(V)$ with non-equivariant degree
$d$ if and only if $d\equiv \pm 1\, ({\rm mod}\, 3)$. 

{\rm (ii).}
There is no $G$-map $P\to S(V)$.
}
\begin{proof}
(i). For the
non-equivariant degree of a $G$-map is congruent (mod $3$)
to the degree of its restriction $S(\Rr )\to S(\Rr )$ to
the fixed subspace $S(V^H)$ of the cyclic subgroup $H$
of order $3$.  And we can realize a map of degree $\pm (3r+1)$ 
by the join of $t\mapsto \pm t: S(\Rr )\to S(\Rr )$ and
$z\mapsto z^{3r+1}: S(\Cc )\to S(\Cc )$.
\par\noindent
(ii). Let $\alpha$ and $\beta$ be the real $1$- and $2$-dimensional
vector  bundles
over $X=P/G=S(\Rr^2\otimes A)/\langle g^3\rangle \times
S(\Rr^3\otimes B)/\langle g^2\rangle$ 
associated with the modules $A$ and $B$. So
$\xi =\alpha\oplus\beta$.

Now, writing $\tilde \Zz$ for integer coefficients twisted by the
Hopf line bundle over $S^1$, we have,
for any compact ENR $Y$, an equivalence
$$
H^i(S^1 \times Y;\,\tilde\Zz) =H^{i-1}(Y;\,\Ff_2),\quad
i\in\Zz ,
$$
given by reduction (mod $2$) and the projection to the second factor
$$
H^i(S^1 \times Y;\,\tilde\Zz) \to H^i(S^1;\,\Ff_2)
=H^i(Y;\,\Ff_2)\oplus H^{i-1}(Y;\,\Ff_2)\to H^{i-1}(Y;\,\Ff_2).
$$
(This is the classical statement that
the real projective plane is a mod $2$ Moore space.
Given the naturality, it is enough to check when $Y$ is a point.)

There is only one non-orientable $2$-dimensional real
vector bundle over $X$, namely $\Rr\oplus\alpha$, because
the group $H^2(X;\,\Zz (\alpha ))=
H^1(S(\Rr^3\otimes B)/\langle g^2\rangle ;\,\Ff_2)$ 
that classifies such bundles is zero. Since the Pontryagin class
$p_1(\beta )\in H^4(X;\,\Zz )$ is non-zero, 
$\xi$ is not isomorphic to $\Rr\oplus (\Rr\oplus\alpha )$. 
\end{proof}
\Thm{Suppose that for each prime $p$ dividing $\# G$ the fixed subspace
$V^H$ of a Sylow $p$-subgroup $H$ is non-zero and, if $n$ is odd, 
of dimension greater than $1$.
Then there exists a $G$-map $S(V)\to S(V)$ with non-equivariant degree
equal to $0$.
}
\begin{proof}
Indeed, since the definition of $\rho_r$, Definition \ref{rho},
is equivariant with respect to the action of the orthogonal group
of $V$,
for any even-dimensional $G$-module $V$ there is
(as was already shown in \cite{lyle}) a non-negative
integer $N_G(V)$ such that there is a $G$-map $S(V)\to S(V)$
with non-equivariant degree $1-d$ if and only if $d$ is divisible
by $N_G(V)$. And then it follows from the equivariant version of
Proposition \ref{mvs} that, if $V^H\not=0$ for a Sylow $p$-subgroup
$H$, the prime $p$ does not divide $N_G(V)$. 
Hence, if this holds for all primes, we must have $N_G(V)=1$.

The proof when $n$ is odd uses the equivariant version of Proposition \ref{mvs1}.
\end{proof}
\Prop{\label{pgrp}
Suppose that $G$ is a $p$-group and that $V^G=0$. Then any $G$-map 
$S(V)\to S(V)$ has non-equivariant degree 
congruent to $1\, ({\rm mod}\, p)$.

If $n$ is even, then $N_G(V)$ is divisible by $p$.
If $n$ is odd, $p=2$.
}
\begin{proof}
In general, for the non-trivial $p$-group $G$,
the degree of a $G$-map $f : S(V)\to S(V)$ is
congruent $({\rm mod}\, p)$ to the degree of its restriction
$f^G : S(V^G)\to S(V^G)$ to the fixed subspace.
(See, for example, \cite{tD}. The map $f$ has an equivariant
degree in the Burnside ring $A(G)$ of $G$, the Grothendieck
group of finite $G$-sets. And so the assertion reduces to
the fact that the cardinality of a finite $G$-set is
congruent $({\rm mod}\, p)$ to the cardinality of the
fixed subset.) 

If $n$ is odd, the map $-1 : S(V)\to S(V)$ has degree $-1$.
\end{proof}
\Prop{\label{pfrp2}
Suppose that $G$ is a group of order $2p^s$, where $p$ is an odd prime,
$s\geq 1$, and that $V^G=0$. Suppose further that for a Sylow
$p$-subgroup $H$, $\dim V^H =1$.  
Then any $G$-map $S(V)\to S(V)$ has non-equivariant degree 
congruent to $\pm 1\, ({\rm mod}\, p)$.

If $n$ is even, then $N_G(V)$ is divisible by $p$.
}
\begin{proof}
Observe that $H$ is a normal subgroup and the quotient
group $G/H$ acts on $V^H$ as $\pm 1$.
The degree of a $G$-map $f : S(V)\to S(V)$ is congruent
$({\rm mod}\, p)$ to the degree of its restriction
$f^H : S(V^H)\to S(V^H)$. Since $f^H$ commutes with the
action of $G/H$, we must have $f^H=\pm 1$.

When $n$ is even, for each $r\in\Zz$ there is a $G$-map
$S(V)\to S(V)$ with degree $1-rN_G(V)$, and so the product 
$rN_G(V)(rN_G(V)-2)$ must be divisible by $p$. 
This is only possible if $p$ divides $N_G(V)$.
\end{proof}
\Sect{\label{sym}
Symmetric groups}
For a subset $J\subseteq I=\{1,\ldots ,q\}$,
as in the Intrduction,
we write $e_J=\sum_{i\in J} e_i\in\Rr [I]$ and $\bar e_J\in L$ for its
image in the quotient $L=\Rr [I]/\Rr e_I$.
Consider a prime $p\leq q$ and write $q=\sum_ra_rp^r$, 
where $0\leq a_r<p$.
A Sylow $p$-subgroup $H$ of $\SS_q$ determines a
partition 
$$
I =\bigsqcup_{i=1}^{\alpha_p(q)} J_i
$$ 
of $I$ into $\alpha_p(q)=\sum_r a_r$
subsets each of order a power of $p$
with $a_r$ of order $p^r$.
The fixed subspace $\Rr [I]^H$ has dimension $\alpha_p(q)$ with
the basis $e_{J_i}$. So $L^H$ has dimension $\alpha_p(q)-1$ and is 
spanned by the $\bar e_{J_i}$ subject to the single relation
$\sum_i \bar e_{J_i}=0$.
The order of a Sylow $p$-subgroup is the highest power, 
$p^{\nu_p(q!)}$, of $p$ that divides $\#\SS_q$, 
where the $p$-adic valuation
$\nu_p(q!)$ of $q!$ is equal to $(q-\alpha_p(q))/(p-1)$.

As in the statement of
Corollary \ref{main}, we consider a principal $\SS_q$-bundle
$P\to X$ over a compact ENR and write $\lambda$ for the 
$(q-1)$-dimensional vector bundle $P\times_{\SS_q} L$ over $X$.
Let $\mu$ be some odd-dimensional Euclidean vector bundle over $X$.
\Lem{\label{D}
Suppose that $p$ is a prime and that $q$ is
not a power of $p$.
If $q$ is odd, then $p$ does not divide $N(\lambda )$;
if $q$ is even, then $p$ does not divide $N(\lambda\oplus\mu )$.
}
\begin{proof}
Choose a point $v\in S(L)$ that is fixed by some Sylow $p$-subgroup
or, equivalently, so that its $\SS_q$-orbit $D$ has size $d$ which is
not divisible by $p$.

For example,
if $p^s$ is the highest power of $p$ dividing $q$,
we could take $v=\bar e_J/\| \bar e_J\|$, where $J\subseteq I$
has cardinality $p^s$.
Then $D=\{ \bar e_K/\| \bar e_K\|  \st K\subseteq I,\, \# K=p^s\}$
and $d=\binom{q}{p^s}$.

Take $\tilde X =P\times_{\SS_q} D
\subseteq P\times_{\SS_q} S(L) =S(\lambda )
\subseteq S(\lambda \oplus\mu)$.
By Proposition \ref{mvs}, $S(\lambda )$ and $S(\lambda\oplus\mu )$
admit self-maps of degree $1-d$.
So, if $q$ is odd, $N(\lambda )$ divides $d$, and, if $q$ is even,
$N(\lambda\oplus\mu )$ divides $d$.
\end{proof}
\Lem{\label{Np}
Suppose that $q$ is a power $q=p^s$, $s\geq 1$, of a prime $p$.
If $p$ is odd, then $N(\lambda )$ is either $1$ or $p$;
if $p=2$, then $N(\lambda\oplus\mu )$ is either $1$ or $2$.
}
\begin{proof}
For we can take $D=\{ \bar e_J\st \# J=p^{s-1}\}$.
Then $\# D=\binom{p^s}{p^{s-1}}$ is divisible by $p$ but not $p^2$.
(For $\alpha_p(p^{s-1})+\alpha_p(p^s-p^{s-1})-\alpha_p(p^s)=p-1$.)
\end{proof}
\Prop{\label{ex1}
Suppose that $q=p^s$ is a prime power, $s\geq 1$, and $k\geq 1$.
Then there exists a principal $\SS_q$-bundle $P\to X$
admitting no $\SS_q$-map $P\to S(\Rr^k\otimes L)$.
}
\begin{proof}
Let $G$ be an elementary abelian $p$-group of order $p^s$.
Fixing a bijection $I\to G$, we can identify $G$ with
a subgroup of $\SS_q$ such that the restriction $V$ of
$\Rr^k\otimes L$ to $G$ satisfies $V^G=0$.
According to Proposition \ref{dickson} there exists
a principal $G$-bundle $Q$ such that there is no $G$-map
$Q\to S(V)$. Then there is no $\SS_q$-map $P=Q\times_G\SS_q
\to S(L)$.
\end{proof}
A similar argument to that in Proposition \ref{counter}
shows that in Corollary \ref{main} when
$k=1$ and $q$ is even some additional condition is essential.
\Prop{\label{ex2}
Take $q=2p^s$, where $p>2$ is an odd prime and $s\geq 1$.
Then there exists a principal $\SS_{2p^s}$-bundle $P\to X$
such that there is no $\SS_{2p^s}$-map $P\to S(L)$.
}
\begin{proof}
Let $G=C_2\times H$ be the product of a cyclic group of order $2$
and an elementary abelian $p$-group $H$ of order $p^s$. 
Fix a bijection $I\to G$ to identify $G$ with a subgroup of 
$\SS_{2p^s}$
and the restriction of $L$ to $G$ with
$$
V=A\oplus B\oplus (A\otimes B)
$$
where $A=\Rr$ is the non-trivial $C_2$-module with the generator
acting as $-1$ and $B=\Rr [H]/\Rr$, of dimension $p^s-1$,
is the reduced regular representation of $H$.

Taking $H=C_p\times\cdots\times C_p$ (where again $C_p$ is
the group of $p$th roots of unity in $\Cc$), let
$R$ be the free $H$-space $S(\Cc^m)\times\cdots\times S(\Cc^m)$
with the action $(h_i)\cdot (v_i)=(h_iv_i)$ for $m$ sufficiently 
large, say $m=2p^s-1$.
Let $Q$ be the free $G$-space 
$S(\Rr^2\otimes A)\times R$ of dimension $(2m-1)s+1$,
and let $\alpha$ and $\beta$ be the real line bundle and
$(p^s-1)$-dimensional real vector bundle over 
$X=Q/G= S(\Rr^2\otimes A)/C_2 \times R/H$
associated with the modules $A$ and $B$ respectively.
We take $P=Q\times_G\SS_{2p^s}$, so that over $X=P/\SS_{2p^s}$
$$
\lambda =\alpha\oplus \beta \oplus (\alpha\otimes\beta ).
$$

Suppose that $\eta$ is a $(2p^s-2)$-dimensional non-orientable real
vector bundle over $X=Q/G=S^1\times R/H$. 
Now, as in the proof of Proposition \ref{counter},
the group 
$$
H^{2(p^s-1)}(X;\,\Zz (\alpha ))=H^{2p^s-3}(R/H ;\,\Ff_2)
$$
is trivial.
So $e(\eta )\in H^{2(p^s-1)}(X;\,\Zz (\alpha ))$
is zero. 
Recall that, for a real vector bundle $\xi$ of even dimension $2n$,
the $n$th Pontryagin class is equal to the square of the Euler class:
$p_n(\xi )=(-1)^nc_{2n}(\Cc\otimes \xi)=(-1)^ne(\Cc\otimes\xi)
=e(\xi\oplus\xi )=e(\xi )^2$.
(The sign relates the orientations on $\Cc\otimes\xi$ 
and $\xi\oplus\xi$.)
So $p_{p^s-1}(\eta )=e(\eta )^2\in H^{4(p^s-1)}(X;\,\Zz )$ 
is zero.
But $p_{p^s-1}(\lambda )$ is non-zero, 
because 
$$
p_{p^s-1}(\beta\oplus\beta )=e(\beta\oplus\beta )^2=
e(\beta )^4\not=0,
$$
by the choice of $m$.
So $\lambda$ cannot be isomorphic to $\Rr\oplus\eta$.

There is, thus, no $G$-map $Q\to S(V)$ and no
$\SS_{2p^s}$-map $P\to S(L)$.
\end{proof}
But it is shown in \cite{AK1} that these (Propositions \ref{ex1} and
\ref{ex2}) are the only exceptional 
cases.
\Thm{{\rm (\cite[Theorem 1.1]{AK1}).}
Suppose that $q$ is even, not a power of $2$,
and not twice a power of an odd prime.
Then, for any principal $\SS_q$-bundle $P\to X$,
there exists an $\SS_q$-equivariant map $P\to S(L)$.
}
The first example is $q=12=2^3+2^2=3^2+3^1$.
\begin{proof}
Let $\Gamma$ denote the set of primes $p\leq q$
and $\Gamma_2\subseteq \Gamma$ the set of primes $p\leq q$ with
$\alpha_p(q)=2$.

For each prime $p\in\Gamma$, choose a partition $I=\bigsqcup J_i$
as above into $\alpha_p(q)$ subsets, and let $G_p$ be the
subgroup of $G=\SS_q$ that preserves each subset $J_i$:
$G_p =\prod_i \SS_{J_i}\leq\SS_I$. 

If $\alpha_p(q)>2$, we can choose two points $x_p$ and 
$y_p\not=\pm x_p$ in
$S(L)$, of the form $\sum_i t_i\bar e_{J_i}$, where $t_i\not=t_j$ for
$i\not=j$ and $\sum_i t_i=0$, with stabilizer exactly $G_p$.
Since $G_p$ contains a Sylow $p$-subgroup of $\SS_q$, the 
the order $d_p$ of the $\SS_q$-orbit of $x_p$ is not divisible by $p$.

If $\alpha_p(q)=2$, say $q=p^s+p^t$,
where $s>t$ we choose a single point $x_p$ with stabilizer
exactly $G_p$. Note that, in this case, 
$d_p=\# (G/G_p)=\binom{p^s+p^t}{p^t}$ is congruent to $1$ (mod $p$).
(For $(1+T)^{p^s+p^t}=(1+T^{p^t})^{p^{s-t}+1}$ in $\Ff_p [T]$.)
Also, $p$ divides $d_{p'}$ for $p'\not=p$.
(If $r_1+\ldots +r_a=p^s+p^t$, where $r_i\geq 1$,
then $\nu_p((r_1+\ldots +r_a)!/(r_1!\cdots r_a!))=
(\alpha_p(r_1)+\ldots +\alpha_p(r_a)-2)/(p-1)$.)

Notice that the $G$-orbits of the points $x_p$ are disjoint.
We write $\tilde X_p = P\times_G Gx_p$ and $\tilde X=\bigsqcup_p
\tilde X_p\subseteq S(\lambda )$.
Following the notation of Proposition \ref{mvs1},
we have a vector bundle $\eta$ over $\tilde X$ of even dimension
$q-2$.
Let $\eta_p$ be the restriction of $\eta$ to $\tilde X_p$.
Then $N(\eta_p)=1$ if $p\notin\Gamma_2$, because $y_p$ gives a
section of $S(\eta_p)$, and $N(\eta_p)=1$ or $p$ if $p\in\Gamma_2$,
by Lemma \ref{Np}.

Now, (following the argument in \cite{AK1}) we notice that
$c=1-\sum_{p\in\Gamma}d_p$ is divisible by each $p\in\Gamma_2$.
We can write $1=\sum_{p\in\Gamma} n_pd_p$, and
so $1=\sum_{p\in\Gamma}(1+cn_p)d_p =
\sum_{p\in\Gamma}m_pd_p$, where $m_p=1+cn_p$ is congruent
to $1 \, ({\rm mod}\, p)$ if $p\in\Gamma_2$.

The proof is completed by applying Proposition \ref{mvs1}.
\end{proof}
\Sect{\label{stable}
Euler classes and stable homotopy theory}
Stable cohomotopy will be denoted by $\omega^*$, and we make the abbreviation 
$$
\omega^0(X;\,-\xi )=\tilde\omega^0(X^{-\xi })
$$ 
for the reduced stable cohomotopy of the Thom space $X^{-\xi}$ of the virtual bundle $-\xi$. Thus $\omega^0(X;\, -\xi )$ may be identified
with the stable cohomotopy group
$\omega^{N-n}(X^{\xi^\perp})$ where $\xi^\perp$ is the 
orthogonal complement of $\xi$ included in a trivial bundle
$X\times\Rr^N$ of sufficiently high dimension $N$.
Similar notation in cohomology $H^*$ with integral
coefficients would write $H^0(X;\,-\xi )$, using the Thom
isomorphism, for $H^n(X;\,\Zz (\xi ))$. 
The stable cohomotopy Euler class $\gamma (\xi )\in\omega^0(X;\, -\xi)$ of $\xi$ is defined like the classical Euler class
$e(\xi )\in H^n(X;\, \Zz (\xi ))$. Further details may be found,
for example, in \cite[II, Section 4]{fht}.

The stable cohomotopy Euler class 
$$
\gamma (\xi )\in\omega^0(X;\,-\xi ),
$$
like $e(\xi )$,
is zero if the sphere bundle $S(\xi )$ has a section.
In a metastable range the converse is true, by Freundenthal's suspension theorem.
\Prop{\label{freuden}
Suppose that $\dim X < 2(\dim\xi -1)$. If $\gamma (\xi )=0$,
then $S(\xi )$ admits a section.
\qed
}
The vanishing of the cohomology Euler class
$e(\xi )\in H^n(X;\,\Zz (\xi ))$ is sufficient for the existence of
a section if $\dim X \leq \dim \xi$.
\Prop{\label{lift}
{\rm (See, for example, \cite[II, Lemma 12.45]{fht}).}
Suppose that $\pi : \tilde X \to X$ is a $d$-fold cover.
If the pullback $S(\pi^*\xi )\to \tilde X$ of the sphere bundle
of $\xi$ admits a section, then
$$
d^k\gamma (\xi ) =0\in \omega^0(X;\, -\xi )
$$
for $k\geq 1$ sufficiently large.
}
\begin{proof}
Since $S(\pi^*\xi )$ admits a section,
the stable cohomotopy Euler class 
$\pi^*\gamma (\xi )=
\gamma (\pi^*\xi)\in \omega^0(\tilde X;\, -\pi^*\xi )$ is zero.
Hence $\pi_!\pi^*\gamma (\xi )= \pi_!(1)\,\gamma (\xi )\in
\omega^0(X;\,-\xi )$ is zero, where $\pi_!:\omega^0(\tilde X)\to
\omega^0(X)$ is the transfer (or direct image) homomorphism.
Now $\pi_!(1)\in\omega^0(X)$ can be written as $d-x$ where
$x$ is nilpotent.  Thus $d\gamma (\xi )=x\cdot \gamma (\xi)$,
and $d^k\gamma (\xi )=x^k\cdot\gamma (\xi )$ is zero
for large $k$.
\end{proof}
\Rem{The same argument in cohomology rather than stable cohomotopy
shows that $d\cdot e(\xi )=0\in H^n(X;\,\Zz (\xi ))$.
}
\Rem{If $\tilde X\to X$ is embedded in $S(\xi )$, so defining
a $d$-valued section, then $S(\pi^*\xi )$, as we observed in Section 
\ref{noakes}, clearly has a section.
}
\Prop{Suppose that $S(\xi )$ has a self-map $f$ of degree $1-d$.
Then 
$$
d^k\gamma (\xi ) =0\in \omega^0(X;\, -\xi )
$$
for $k\geq 1$ sufficiently large.
}
\begin{proof}
We have $\gamma (\xi )=[f]\cdot \gamma (\xi )$,
where $[f]\in\omega^0(X)$ is the class represented by $f$.
Now $[f]=1-d+x$, where $x$ is nilpotent. So the result
follows as before.
\end{proof}
\Cor{\label{odd}
If $\dim\xi$ is odd, then $2^k\gamma (\xi )=0$ for
$k$ large.
}
\begin{proof}
The antipodal map $-1: S(\xi )\to S(\xi )$ has degree $-1$.
\end{proof}
\Prop{\label{mvs2}
If $f$ is the map constructed in Propositions \ref{mvs}
and \ref{mvs1}
from a $d$-valued section of $S(\xi )$ specified by
a $d$-fold cover $\pi :\tilde X\to X$ and fibrewise embedding
$\iota : \tilde X \into S(\xi )$ over $X$
and a map $h: S(\eta )\to S(\eta )$, then 
$$
[f]=1-\pi_!([h])\in \omega^0(X),
$$ 
where $[h]\in\omega^0(\tilde X)$ is the class represented by $h$.
}
\begin{proof}
The section $s_{\tilde X,h}$ in Construction \ref{thom1} realizes 
the fibrewise Umkehr homomorphism 
$$
\iota_! : \omega^0(\tilde X)\to \omega^0(S(\xi );\,-\tau_XS(\xi )),
$$
where $\tau_XS(\xi)$ is the fibrewise tangent bundle of
$S(\xi )\to X$. (See, for example, \cite[II, Chapter 3]{fht}.)
Thus, $[s_{\tilde x,h}]=\iota_!([h])$.
The projection $\pi$ is the composition of $\iota$ and
the projection $\sigma : S(\xi )\to X$.
Hence, $\sigma_![s_{\tilde X,h}]=\pi_![h]$.

Recall that two sections $s_0,\, s_1$ of the
pullback $S(\sigma^*\xi )$ determine
a difference class
$\delta (s_0,s_1)\in\omega^{-1}(S(\xi );\, -\xi)$
constructed as an obstruction to the existence of a
homotopy between the two sections.

The class 
$$
[s_{\tilde x,h}]\in\omega^0(S(\xi );\, -\tau_X S(\xi ))=
\omega^{-1}(S(\xi );\, -\xi )
$$
can be expressed as $\delta (s_\emptyset ,s_{\tilde X,h})$,
where $s_\emptyset$ is the section at $\infty$ of
$(\tau_XS(\xi ))^+_{S(\xi )} =S(\Rr\oplus \tau_XS(\xi ))$.
In the Gysin sequence
$$
\cdots\to
\omega^{-1}(S(\xi );\,-\xi )\Rarr{\delta}{}
\omega^0(X)\Rarr{\gamma (\xi )\cdot}{}
\omega^0(X;\,-\xi )=\omega^0(D(\xi ),S(\xi );\, -\xi)
\to\cdots
$$
of the sphere bundle $S(\xi )$, that is,
the exact sequence of the pair $(D(\xi ),S(\xi ))$,
the coboundary $\delta$ coincides with the Umkehr
map $\sigma_!$ and so maps the difference class to
$\sigma_![s_{\tilde X,h}]=\pi_![h]$.

Given any two self-maps $f_0,\, f_1 :S(\xi )\to S(\xi )$,
the difference class $\delta (s_0,s_1)$  of the corresponding
sections $s_0$ and $s_1$ of $S(\sigma^*\xi )$ maps,
in the Gysin sequence, to 
the difference $\gamma (\xi ,s_0)-\gamma (\xi ,s_1)$
of the relative Euler classes
$\gamma (\xi ,s_0)$ and $\gamma (\xi ,s_1)$
in $\omega^0(D(\xi ),S(\xi );\, -\xi )$.
(See \cite[II, Lemma 4.12]{fht}.)
The relative Euler class $\gamma (\xi ,s_i)$ is an obstruction
to extending $f_i : S(\xi )\to S(\xi )$ to a map
$D(\xi )\to S(\xi )$ and corresponds to $[f_i]\in
\omega^0(X)$. 
Applying this to the self-maps $f_0=1$ and $f_1=f$, which correspond,
respectively, to the sections $s_\emptyset$ and $s_{\tilde X,h}$,
we obtain the required identity $1-[f]=\pi_![h]$.
\end{proof}

Thus far, the arguments have applied to a general vector bundle $\xi$ over a compact ENR $X$. We now specialize to the case that
$\xi =P\times_G V$ as in the Introduction.
\Cor{\label{prime}
Let $p$ be a prime. Suppose that for a Sylow $p$-subgroup 
$H\leq G$ the fixed subspace $V^H$ is non-zero.
Then the $p$-local Euler class $\gamma (\xi)\in
\omega^0(X;\, -\xi )_{p)}$ is zero.
}
\begin{proof}
Take $\pi :\tilde X \to X$ in Proposition \ref{lift}
to be the covering projection $P/H\to P/G$
with order $d=\# (G/H)$ prime to $p$.
\end{proof}
\Cor{\label{global}
Suppose that for each prime $p$ dividing $\# G$ the fixed subspace
$V^H$ of a Sylow $p$-subgroup $H$ of $G$ is non-zero. Then
the stable cohomotopy Euler class 
$\gamma (\xi )\in\omega^0(X;\, -\xi )$ is zero.
}
\begin{proof}
Since, by Corollary \ref{prime},
$\gamma (\xi )$ vanishes at any prime $p$, we conclude that the stable
cohomotopy Euler class of $\xi$ is zero.
\end{proof}
\Cor{\label{euler}
Suppose that for each prime $p$ dividing $\# G$ the fixed subspace
$V^H$ of a Sylow $p$-subgroup $H$ of $G$ is non-zero. 
If $\dim X < 2(n-1)$,
then the sphere bundle $S(\xi  )\to X$ has a section.
} 
\begin{proof}
Since $\gamma (\xi )=0\in \omega^0(X;\, -\xi )$,
by Corollary \ref{global},
the conclusion follows from Proposition \ref{freuden}.
\end{proof}
\begin{proof}[Proof of Theorem \ref{gen}(b)]
If $n$ is odd, $\gamma (\xi )$ is $2$-torsion,
by Corollary \ref{odd}.
So it enough to show that the $2$-local Euler class in
$\omega^0(X;\, -\xi )_{(2)}$ is zero.
\end{proof}
\Thm{Let $r\geq 1$ be a positive integer such that
$\dim X < 2(n-r)$ and $r<n$. 
Suppose that, for each prime $p$ dividing $\# G$ if $n$ is even,
for $p=2$ if $n$ is odd and $\# G$ is even,
the dimension of the fixed subspace $V^H$ of a Sylow $p$-subgroup
$H$ of $G$ is at least $r$. Then there exists a $G$-map
$$
P \to \O (\Rr^r,V)
$$
to the Stiefel manifold of orthogonal $r$-frames in $V$.
}
\begin{proof}
In this range of dimensions the Stiefel bundle 
$$
\O (\Rr^r,\xi )=P\times_G \O (\Rr^r,V)\to X
$$
has a section 
if and only if the stable cohomotopy Euler class
of the tensor product $\eta\otimes\xi$ 
of the Hopf line bundle $\eta$ over the real
projective space $P(\Rr^r)$ of $\Rr^r$ with $\xi$,
$$
\gamma (\eta\otimes\xi )\in 
\omega^0(P(\Rr^r)\times X;\, -\eta\otimes\xi),
$$
is zero. (See, for example, \cite[Theorem 1.5]{KO}.)

For a Sylow $p$-subgroup $H$ of $G$, consider the finite
cover $\pi : \tilde X =P/H\to X=P/G$ with order prime to $p$.
The pullback $\pi^*\xi$ has a trivial subbundle of
dimension $r$ and $S(\eta\otimes\pi^*\xi)$ has a section. 
So $\gamma (\eta\otimes\pi^*\xi )
\in\omega^0(P(\Rr^r)\times\tilde X;\, -\eta\otimes\pi^*\xi )$
is zero. From Proposition \ref{lift} applied to $\eta\otimes\xi$,  we
deduce that $\gamma (\eta\otimes \xi )$ is zero at the prime
$(p)$.
Hence, using Corollary \ref{odd} if $n$ is odd, we conclude that
$\gamma (\eta\otimes\xi )$ is zero.
\end{proof}

\end{document}